\documentclass{ifacconf}

\usepackage{graphicx}      
\usepackage{natbib}        

\usepackage{amsmath,amssymb,mathtools}
\usepackage{booktabs}
\usepackage{multirow}

\usepackage{color}
\definecolor{forestgreen}{rgb}{0.33,0.61,0.34}

\DeclareMathOperator*{\minimize}{minimize}

\begin{document}
\begin{frontmatter}

\title{Static Output Feedback Synthesis of Time-Delay Linear Systems via Deep Unfolding\thanksref{footnoteinfo}} 

\thanks[footnoteinfo]{This work was supported by JSPS KAKENHI Grant Numbers JP21H01353 and JP22H00514.}

\author[OU,NAIST]{Masaki Ogura} 
\author[NAIST]{Koki Kobayashi} 
\author[NANZAN,NAIST]{Kenji Sugimoto} 

\address[OU]{Osaka University, Suita, Osaka 565-0871, Japan (e-mail: m-ogura@ist.osaka-u.ac.jp).}
\address[NAIST]{Nara Institute of Science and Technology, Ikoma, Nara 630-0192, Japan}
\address[NANZAN]{Nanzan University, Nagoya, Aichi 466-8673, Japan}

\begin{abstract}                
We propose a deep unfolding-based approach for stabilization of time-delay linear systems. Deep unfolding is an emerging framework for design and improvement of iterative algorithms and attracting significant attentions in signal processing. In this paper, we propose an algorithm to design a static output feedback gain for stabilizing time-delay linear systems via deep unfolding. Within the algorithm, the learning part is driven by NeuralODE developed in the community of machine learning, while the gain verification is performed with linear matrix inequalities developed in the systems and control theory. The effectiveness of the proposed algorithm is illustrated with numerical simulations. 
\end{abstract}


\end{frontmatter}

\vspace{8mm}\section{Introduction}\vspace{-2mm}

Stabilization is one of the fundamental issues in the context of time-delay systems~\citep{Fridman2014b} and has been actively investigated in the literature. For example, the Smith predictor approach \citep[e.g.,][]{Mirkin2011} allows us to reduce the stabilization problem into a delay-free problem. The eigenvalue-based approach \citep[e.g.,][]{Michiels2005} provides us with intuitive and generally applicable methodologies for stabilization. Recently, \citet{Barreau2018} presented a methodology for designing a static feedback gain for stabilization with iterative linear matrix inequalities derived from Lyapunov-Krasovskii functionals. 

The objective of this paper is to present a \emph{deep unfolding}-based approach for stabilization of time-delay linear systems. Deep unfolding is a learning-based approach for design and improvement of iterative algorithms, and has been successfully employed in signal processing for wireless communications~\citep{Jagannath2021}. Recently, in the context of the systems and control engineering, \citet{Kishida2021b} have demonstrated the effectiveness of  deep unfolding in nonlinear model predictive control. The core idea of  deep unfolding is in regarding an iterative algorithm (or, a dynamical system as a special case) as a signal flow graph, which can then be efficiently trained by using techniques available in the field of machine learning. In this paper, motivated by deep unfolding, we develop an algorithm for finding a stabilizing static output feedback gain for time-delay linear systems, and numerically illustrate the algorithm's effectiveness via numerical simulations.

\vspace{8mm}\section{Problem and algorithm}\vspace{-2mm}

Let $n$, $m$, and $p$ be positive integers. Let $h >0$ be a constant. For matrices~$A \in \mathbb R^{n\times n}$, $B \in \mathbb R^{n\times m}$, and $C\in \mathbb R^{p\times n}$, let us consider the time-delay linear system
\begin{equation}\label{eq:def:sigma}
\Sigma: \begin{cases}
    \dot x(t) = Ax(t) + Bu(t-h), 
    \\
    y(t) = Cx(t)
\end{cases}
\end{equation}
with the initial condition $x\vert_{[-h, 0]} = \phi \in \mathcal C([-h, 0], \mathbb R^n)$, where $\mathcal C([-h, 0], \mathbb R^n)$ denotes the normed space of~$\mathbb R^n$-valued continuous functions defined on~$[-h, 0]$. Consider the static output feedback 
    $u(t)  = Ky(t)$, 
where $K \in \mathbb R^{m\times p}$ is the feedback gain to be designed. We say that the closed loop system, denoted by~$\Sigma_K$, is \emph{asymptotically stable} if $x(t)$ converges to $0$ as $t\to\infty$ for any initial condition~$\phi$. 

Our objective is to find a stabilizing feedback gain $K$ given the coefficient matrices $A$, $B$, and $C$ as well as the delay length~$h$ of~$\Sigma$. 
We start the derivation of our algorithm from the following trivial but important fact. 

\begin{lem}\label{lem:trivial}
For $K \in \mathbb{R}^{m \times p}$ and $t>0$, define
\begin{equation}
    \ell_t(K) = \sup_{\lVert \phi\rVert = 1} \lVert x(t)\rVert, 
\end{equation}
where $x$ denotes the trajectory of~$\Sigma_K$ with initial condition~$\phi$. Then, $\Sigma_K$ is asymptotically stable if and only if $K$ minimizes the function $K\mapsto \limsup_{t\to\infty}\ell_t(K)$ with the minimum value $0$. 
\end{lem}

Motivated by Lemma~\ref{lem:trivial}, we introduce the following relaxed optimization problem: 
\begin{equation}\label{eq:relaxed}
    \minimize_{K\in\mathbb{R}^{m\times p}}\quad E[\lVert x(T) \rVert], 
\end{equation}
where $T > 0$ is a constant, $x$ denotes the trajectory of~$\Sigma_K$ with a  \emph{random} initial condition~$\phi$, and $E[\cdot]$ denotes the mathematical expectation with respect to~$\phi$. From Lemma~\ref{lem:trivial}, if $T$ is large and the support of the distribution spans a large subspace, then we can expect that solving the optimization problem~\eqref{eq:relaxed} leads to a stabilizing gain~$K$. The effectiveness of this argument is numerically confirmed by~\citet{Kishida2021b} in the context of nonlinear model predictive control. 

A major challenge in finding a (sub-)optimal solution of the optimization problem~\eqref{eq:relaxed} stems from the fact that $x$ is the trajectory of a time-delay system. To overcome this issue, we propose using NeuralODE~\citep{Chen2018b}. Roughly speaking, NeuralODE allows us to `unfold' the dynamics of~$\Sigma_K$ and efficiently evaluate the derivative of the terminal cost~$\rVert x(T)\lVert$. Therefore, starting from arbitrary initial (and deterministic) state~$\phi$ and gain~$K$, we can perform a gradient descent to find a gain making the cost~$\rVert x(T)\lVert$ smaller, as commonly done in deep unfolding~\citep{Jagannath2021}. After finishing the gradient descent, we can theoretically check if the resulting gain is stabilizing by using the LMIs presented by \citet[][Theorem 1]{Barreau2018}. 

Let us present a further detail of the proposed algorithm. Let $M$ and $J$ be positive integers. We draw samples $\phi_1$, \dots, $\phi_J\in \mathcal C([-h, 0], \mathbb R^n)$ of the initial state. Starting from an initial gain $K$, we perform gradient descent on the cost $\lVert x(T/M)\rVert$ (not $\lVert x(T)\rVert$) using the $J$ samples with a pre-specified batch size. We then check if the gain $K$ learned stabilizes $\Sigma$ by using \citep[Theorem 1]{Barreau2018}. If not stabilizing, then we use the learned gain as the initial value to perform further gradient descent on the cost~$\lVert x(2T/M)\rVert$. We repeatedly perform this procedure (called an incremental learning in deep unfolding) until we finish minimizing the terminal cost $\lVert x(T)\rVert$. 

\begin{figure}[tb]
    \centering
    \includegraphics[width=.87\linewidth]{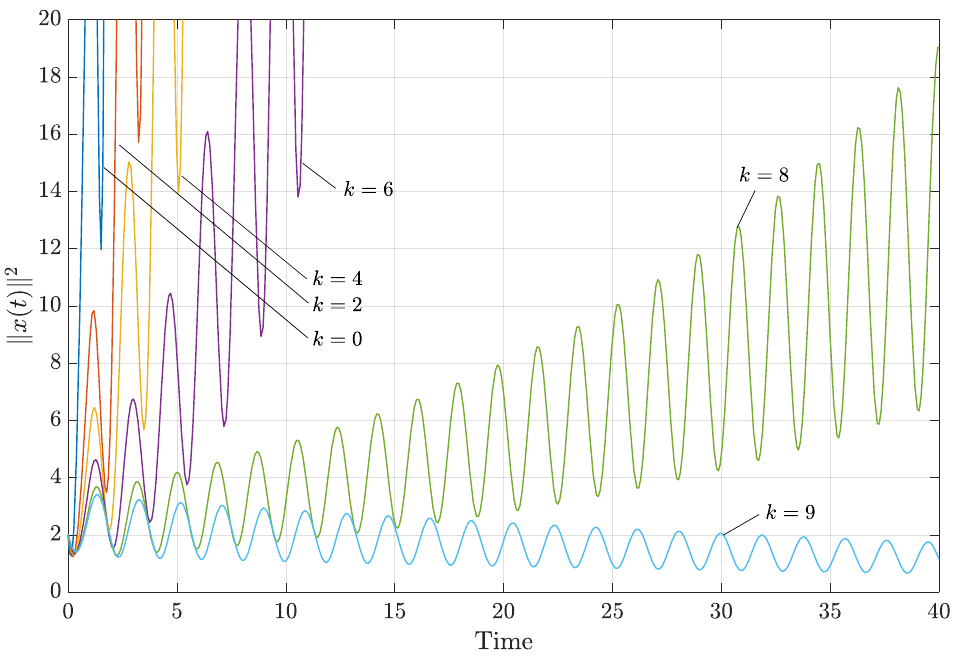}
    \vspace{-2mm}
    \caption{Trajectories of~$\Sigma_K$ with gain $K$ obtained after $k$th step (i.e., after minimizing $\lVert x(k)\rVert$.)}
    \label{fig:my_label}
\end{figure}

\section{Numerical simulations}\vspace{-2mm}

We first consider the time-delay linear system in \citep[Section~VI.B]{Barreau2018}, where the authors suggest the limitations in effectively applying their stabilization algorithm. We set the delay as $h=1$, and we applied the proposed algorithm. We used the parameters $T=20$, $J=10$, and $M=20$. We used \texttt{adam} with learning rate~$0.1$. The proposed algorithm terminated at the $9$th step (i.e., after minimizing $\lVert x((9/M)T)\rVert = \lVert x(9)\rVert$) in less than 4 seconds, with a feedback gain theoretically confirmed to be stabilizing. We illustrate in Fig.~\ref{fig:my_label} how the closed-loop trajectories improved with the incremental learning. 

To quantitatively examine the proposed algorithm's performance, we conduct the following experiment. We consider Scenario~1: $(n, m, p, h) = (4, 1, 2, 0.1)$ and Scenario~2: $(n, m, p, h) = (4, 2, 1, 0.1)$. For each scenario, we randomly generate 100 open-loop systems $\Sigma$ by drawing matrices~$A$, $B$, and $C$ from a distribution. For each system, we apply the proposed algorithm. For comparison, we consider the following two conventional methods; in BMI, we solve the bilinear matrix inequality resulting from~\cite[Theorem~1]{Barreau2018}, while in ILMI, the iterative linear matrix inequality presented in \cite[Section~V]{Barreau2018} is solved. Hence, the proposed algorithm uses the analysis result by~\citet{Barreau2018} for checking stability, while the conventional one use the synthesis algorithms developed by the same authors. For fairness, we use the same hyper-parameter~$N$ within the three algorithms. 

The results of the experiment are summarized in Table~\ref{tbl:}. The proposed algorithm outperforms both the conventional methods in terms of the frequency of stabilization. We remark that the execution time of the proposed algorithm was not necessarily significantly larger than BMI. We also remark that ILMI terminated fast but failed to find a stabilizing feedback gain in any of the scenarios. 

Because the results reported in this paper are preliminary, it is necessary to perform further and thorough comparison between the algorithms. It is also necessary to theoretically analyze properties of the proposed algorithm.

\begin{table}[tb]
\caption{Comparison of algorithms}
\label{tbl:}
\vspace{-1.5mm}
\centering
\includegraphics[width=.85\linewidth]{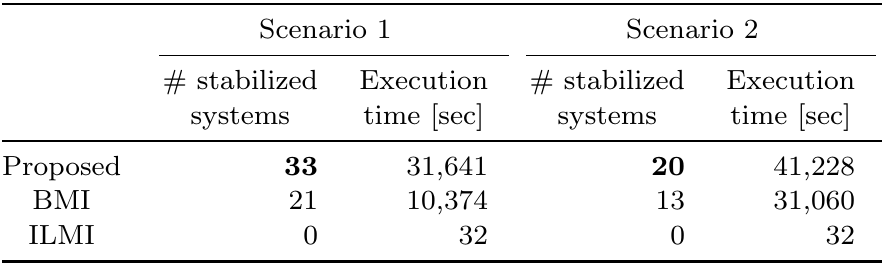}
\end{table}




\end{document}